\documentclass{article}
\usepackage{amssymb,amsmath,latexsym, amsthm}
\usepackage{graphicx}

\newtheorem{theorem}{Theorem}
\newtheorem{lemma}{Lemma}

\newtheorem{proposition}{Proposition}
\newtheorem{definition}{Definition}

\newtheorem{remark}{Remark}

\begin{document}
\title{A non-MRA $C^r$ frame wavelet with rapid decay}

\author{Lawrence Baggett\thanks{The first two named authors were supported by a US-NSF 
Focused Research Group (FRG) grant. }
\\
Department of Mathematics, Campus Box 395\\
 University of Colorado, Boulder, CO 80309-0395\\
\texttt{baggett@euclid.colorado.edu}
\and
Palle Jorgensen\footnotemark[1]
\\
Department of Mathematics, 14 MacLean Hall\\
University of Iowa, Iowa City, IA 52242-1419\\
\texttt{jorgen@math.uiowa.edu}
\and
Kathy Merrill\\
Department of Mathematics, Colorado College\\
 Colorado Springs, CO 80903-3294\\
\texttt{kmerrill@coloradocollege.edu}
\and
Judith Packer\\
Department of Mathematics, Campus Box 395\\
University of Colorado, Boulder, CO 80309-0395\\
\texttt{jpacker@euclid.colorado.edu}}
\date{}
\maketitle

\begin{abstract}
A generalized filter construction is used to build an example of a non-MRA normalized tight frame wavelet
for dilation by 2 in $L^2(\mathbb R)$.  This example has the same multiplicity function as
the Journ\'e wavelet, yet has a $C^{\infty}$ Fourier transform and can be made to be $C^r$ for
any fixed postive integer $r$.
\end{abstract}

\noindent AMS Subject Classification Primary 54C40, 14E20; Secondary 46E25, 20C20

\section{Introduction}
A central problem in the history of wavelet theory has been the search for smooth well-localized
functions whose dilates and translates can be used to obtain all members of a specified function space.
We will use here the simplest definition of an orthonormal wavelet, by taking our context to be dilation by
2 in $L^2(\mathbb{R})$:
\begin{definition}
\label{wavdef}
\indent $\psi\in L^2(\Bbb R)$ is an {\bf orthonormal wavelet} if
 $\{\psi_{j,k}\equiv\sqrt {2} ^j\psi(2^jx-k)\}_{j,k\in\mathbb Z}$ form an orthonormal basis
for $L^2(\mathbb R).$
\end{definition}
Many of the techniques and results used in this paper apply in much wider contexts, to dilation by an expansive matrix
in $L^2(\mathbb{R}^n)$, and to allow for multiple wavelets.  However, for the purposes of creating the example that is the
main thrust of this paper, this more restrictive context will suffice.

The earliest example of an orthonormal wavelet was given by A. Haar in 1910 \cite{Ha}.  Here the wavelet is defined in 
terms of 
characteristic functions of bounded sets,
$$\psi=\chi_{[0,\frac 12)}-\chi_{[\frac 12,1)},$$
and thus is well-localized but not at all smooth.  Another famous classical example,
the Shannon wavelet, has Fourier transform of the wavelet equal to a characteristic
function of a bounded set,
$$\widehat{\psi}=\chi_{[-1,-\frac12)\cup [\frac 12, 1)},$$
and thus has the opposite problem of being smoooth but not well-localized.
In 1988, Ingrid Daubechies \cite{Da2} gave
a construction technique for orthonormal wavelets for dilation by 2 in $L^2(\mathbb{R})$ that are both
smooth and well-localized.  She described how to build, for any integer $r\geq0$ , a wavelet
that is $C^r$ and has compact support.  It can be shown that no wavelet can be both
$C^{\infty}$ and have compact support (see e.g. \cite{HW}).

Daubechies' construction technique used filters in a
multiresolution structure called an MRA.  Most classical wavelets are associated with such
a structure.  J.L. Journ\'e gave the first known example of a wavelet that can not
have an associated MRA (see \cite{Ma}). Journ\'e's wavelet, like the Shannon wavelet, is a {\it wavelet set}
wavelet.  That is, its Fourier transform is a characteristic function.  Since the announcement of the
Journ\'e wavelet, many other non-MRA wavelets have appeared in the literature (see e.g. \cite{DLS}, \cite{BL} \cite{HKL}), most of them also wavelet set
wavelets. These
wavelets are smooth but not well-localized. In fact, it can be shown that
every orthonormal wavelet $\psi$ such that $|\hat{\psi}|$ is continuous and $|\hat{\psi}(\xi)|=O(|\xi|^{-\frac 12 -\epsilon})$
at $\infty$ for some $\epsilon>0$ must be an MRA wavelet (see \cite{HW}).
Thus, in the quest to find wavelets that are both well-localized and smooth, non-MRA wavelets would
seem to be not very useful.  However, if we expand our definition of wavelet slightly, the situation changes.

\begin{definition}\label{frame}
The set of functions $\{\eta_{j}\}$ is a \emph{frame} for
$L^2(\mathbb R)$ if there exist constants $A$ and $B$ such that for each $f\in L^2(\mathbb{R})$ we have
$A\|f\|^2\leq\sum_{j}|\langle f|\eta_{j}\rangle |^2\leq B\|f\|^2$.
The set $\{\eta_{j}\}$ is a \emph{normalized tight frame} or \emph{Parseval frame} if $A=B=1$.
\end{definition}

\begin{definition}\label{framewavelet}
 $\psi\in L^2(\mathbb R)$ is a \emph {(normalized tight) frame wavelet}
or a \emph {Parseval frame wavelet}
for dilation by 2 if $\{\psi_{j,l}\equiv\sqrt {2} ^j\psi(2^jx-l)\}$ form a
normalized tight frame for $L^2(\mathbb R).$
\end{definition}

Note that a normalized tight frame can exhibit redundancy, and
therefore need not be a basis.  Indeed, it can include 0 as one of
its elements.  However, a normalized tight frame $\{\eta_j\}$ does
have the property that every $f\in L^2$ can be recaptured from its
coefficients, $f=\sum\langle f,\eta_j\rangle \eta_j$. (See, e.g. \cite{HW}).
Thus replacing orthonormal
wavelets by frame wavelets does not give up the essential property of
determining all functions in $L^2$ in terms of translates and dilates of the wavelet.

Our main goal in this paper is to produce an example of a non-MRA normalized tight frame
wavelet for dilation by 2 in $L^2(\mathbb R)$ which is both well-localized and smooth, in the sense
that $\psi\in C^r$ (for any fixed positive integer $r$) and $\widehat{\psi}\in C^{\infty}$.  The existence of such a wavelet is somewhat
surprising in view of the result mentioned above for orthonormal wavelets.  While several researchers
have produced examples of frame wavelets that exhibit smoothness and good localization properties (See e.g. \cite{ACM}),
the example in this paper stands out in its relationship to a non-MRA multiresolution structure with a known
multiplicity function.  For such a non-MRA wavelet,
the multiresolution structure will be replaced by a generalized multi-resolution structure called a GMRA,
defined in \cite{HKL}, which we review in Section 2 of this paper.
We note that similar generalizations of multi-resolution structures have also been defined
in \cite{BLi} and \cite{HLPS}.

The Fourier transform is an essential tool in analyzing wavelets; we will use it and a generalization
in the GMRA context to define filters used in wavelet
construction.  We take as the definition of the Fourier transform,
$$\widehat{f}(x)\;=\;\int_{\mathbb R}f(t)e^{-2\pi ixt}dt.$$
Note that on $\widehat{L^2}(\mathbb R)$, translation by $n$ becomes multiplication by $e^{2\pi i n x}$, and dilation
becomes $\widehat{\delta}f(x)=\frac 1{\sqrt2} f(\frac x2).$  We review the filter construction process and its generalization in Section 3 of this paper.
In Section 4 we build our main example.

\section{Generalized multi-resolution analyses}

As mentioned in the introduction, most classical wavelets
have an associated multi-resolution
structure of the following type, first defined by S. Mallat
\cite{Ma}:

\begin{definition}
A \emph{Multiresolution Analysis (MRA)} in $L^2(\mathbb{R})$ is
a collection of closed subspaces $V_j$ that have the following
properties:
\begin{enumerate}
\item $V_j \subset V_{j+1}$
\item $V_{j+1}=\{\delta (f)\equiv \sqrt{2}f(2x)\}_{f\in V_j}$
\item $\cup V_j$ is dense in $L^2(\mathbb R )$ and $\cap V_j = \{0\}$
\item $V_0$ has a \emph{scaling function} $\phi$ whose translates form an orthonormal basis for $V_0$
\end{enumerate}
\end{definition}
However, the Journ\'e wavelet (see \cite{Ma}),
$$\widehat{\psi}=\chi_{[-\frac {16}7,-2)\cup [-\frac 12,-\frac
27)\cup [\frac 27,\frac 12]\cup [2,\frac {16}7)}.$$
has no associated MRA.  This can be shown using the dimension function,
\begin{equation}
\label{dimension}
D(x)\equiv\sum_{k\in\mathbb{Z}}\sum_{j=1}^{\infty}|\widehat{\psi}(2^j(x+k))|^2,
\end{equation}
first introduced by P. Auscher \cite{Au}.
Because an MRA requires the existence of a scaling function, it turns out that a wavelet has an associated MRA 
if and only if the dimension function is identically equal to 1 (see e.g. \cite{HW}).  By removing the scaling
function condition, we
are able to associate a multiresolution structure with the Journ\'e wavelet.
More generally, given any
orthonormal wavelet $\psi\in L^2(\mathbb R)$, the subspaces $\{V_j\}$
defined by $V_j=$ the closed linear span of $\{\psi_{l,k}\}_{l<j}$,
do determine a \textit{generalized}
multiresolution structure, according to the definition below \cite{HKL}:

\begin{definition}
\label{GMRA}
A \emph{Generalized Multiresolution Analysis (GMRA)} is a
collection of closed subspaces $\{V_j\}_{j\in\mathbb Z}$ of
$L^2(\mathbb R)$ such that:
\begin{enumerate}
\item $V_j \subset V_{j+1}$
\item $V_{j+1}=\{\delta(f)\equiv\sqrt{2}f(2x)\}_{f\in
V_j}$
 \item $\cup V_j$ dense in $L^2(\mathbb R)$ and $\cap V_j = \{0\}$
 \item$V_0$ is invariant under translation.
\end{enumerate}
\end{definition}

The definitions of MRA and GMRA differ only in condition (4): An
MRA requires that $V_0$ has a scaling function $\phi$ such that
translates of $\phi$ form an orthonormal basis for $V_0$, while a
GMRA requires only that $V_0$ be invariant under translation by the
integers. In spite of this difference, it is shown in \cite{HKL} that a GMRA has
almost as much structure as an MRA. Translation is a unitary
representation of $\mathbb Z$ on $V_0$, and thus is completely
determined by a multiplicity function $m:[-\frac 12,\frac
12)\mapsto \{0,1,2,\cdots,\infty\}$ describing how many times
each character occurs as a subrepresentation.   The multiplicity function has been shown \cite{W} to be equal to
the dimension function defined in Eq. (\ref{dimension}).
By writing $V_1=V_0\oplus W_0$,
representation theory can be used (see \cite{HKL})
to show that the GMRA has an
associated orthonormal wavelet if and only if the multiplicity
function satisfies a \emph{consistency equation}:
$$m(x)+1=
m(\frac x2)+m(\frac{x+1}2).$$
A GMRA is an MRA iff $m\equiv 1$, and
Journ\'e's famous non-MRA wavelet
example has
$$m(x) = \left\{\begin{array}{ll}
2 & x\in[-\frac 17,\frac 17) \\
1 & x\in\pm[\frac 17,\frac 27)\cup \pm [\frac 37,\frac 12)\\
  0 &
\mbox{otherwise}
\end{array}\right.$$

In every GMRA, there is a unitary equivalence between translation
on $V_0$ and multiplication by exponentials on $\oplus L^2(S_j),$
where $S_j=\{x:m(x)\geq j\}.$  This unitary equivalence plays a
role here similar to that of the Fourier transform in the
classical MRA case. It ensures that a GMRA whose multiplicity function is finite a.e. has
\emph{generalized scaling functions}, $\phi_1,\phi_2,\cdots$ such
that $\{\phi_i(x-l)\}$
form a normalized tight frame for $V_0$.  If $m$ is bounded, with maximum value $c$,
then $c$ generalized scaling functions are required.  If the GMRA is associated with a wavelet set
wavelet, these functions can be chosen
of the form $\widehat{\phi_i}=\chi_{E_i}$, where $\chi_{E_i}$ is the characteristic function
of the set $E_i$, chosen to be congruent mod 1 to $S_i$, and so that $\cup E_i\subset 2\left(\cup E_i\right)$.
(For more details see Theorem 3.5 of \cite{BM}.)  For example, the Journ\'e wavelet has generalized scaling functions
$\widehat{\phi_1}=\chi_{[\frac{-2}7,\frac27)\cup\pm[\frac12,\frac47)}$
and $\widehat{\phi_2}=\chi_{\pm[1,\frac87)}$

\section{Filter constructions}
We begin by recalling the classical filter construction developed by Mallat \cite{Ma} and Meyer
\cite{Me}, and used by Daubechies \cite{Dau} to build $C^r$
wavelets with compact support.
Thus, suppose we have a single wavelet $\psi$ for dilation by 2 in $L^2(\mathbb R),$ with an associated MRA,
and so a scaling function $\phi$ whose translates form an orthonormal basis for $V_0$.

Because $\widehat{\phi}$ and $\widehat {\psi}$ are in $\widehat{V}_{1}$, we can write
$\widehat{\phi}$ and
$\widehat{\psi}$ in terms of exponentials times
the dilate of $\widehat{\phi}$.  That is, there must exist
periodic functions (with period 1)  $h$ and $g$ such that
\begin{equation}
\widehat{\phi}(x)=\frac 1{\sqrt{2}}
h(\frac{x}2)\widehat{\phi}(\frac{x}2)\label{lp}
\end{equation}
 and
 \begin{equation}
\widehat{\psi}(x)=\frac 1{\sqrt{2}}
g(\frac{x}2)\widehat{\phi}(\frac{x}2)\label{hp}.
\end{equation}

The functions $h$ and $g$ are called {\bf low and high pass
filters}.  The roots of this name can be seen by noting that
for the Shannon wavelet, $h=\sqrt 2 \chi_{[-\frac 14,\frac 14)}$
and $g=\sqrt 2 \chi_{[-\frac 12,-\frac 14)\cup[\frac 14,\frac 12)}$ do indeed act by filtering out all but low (for $h$) or
high (for $g$) frequencies.  Because of the orthonormality
conditions satisfied by translates of $\phi$ and $\psi$, all
filters defined by Eqs.(\ref{lp}) and (\ref{hp}) must satisfy
orthonormality-like conditions:
\begin{equation}
|h(x)|^2+|h(x+\frac 12)|^2=2\label{orth1}
\end{equation}
\begin{equation}
|g(x)|^2+|g(x+\frac 12)|^2=2\label{orth2}
\end{equation}
and
\begin{equation}
h(x)\overline{g(x)}+h(x+\frac 12)\overline{g(x+\frac
12)}=0.\label{orth3}
\end{equation}

The classical filter techniques for building wavelets reverse this process
of finding filters from wavelets.  First note that
if $h$ is any periodic function that satisfies Eq.(\ref{orth1}), the function
\begin{equation}\label{classicalg}
g(x)=e^{2\pi ix}\overline{h(x+\frac 12)},
\end{equation}
and $h$ together satisfy all three orthonormality conditions (\ref{orth1}),(\ref{orth2}) and
(\ref{orth3}). (Other choices for $g$
are possible.)  The following classical theorem gives conditions
under which the Fourier transform of a scaling function $\widehat{\phi}$ can then be built by iterating 
Eq.(\ref{lp}), so that the Fourier transform of a wavelet, $\widehat{\psi},$ can be defined using Eq.(\ref{hp}).

\begin{theorem} \label{classical} Let $h$ and $g$ be $C^1$
 functions that satisfy the orthonormality conditions
(\ref{orth1}), (\ref{orth2}), and (\ref{orth3}). Suppose, in
addition, that $h$ is nonvanishing on  $[-\frac 14,\frac 14)$, and
$|h(0)|=\sqrt 2$.  Then:
$$\widehat{\phi}(x)=\prod_{j=1}^{\infty} \frac 1{\sqrt 2}h(2^{-j}x)$$
is a scaling function for an MRA,  and $$\widehat{\psi}(x)=\frac
1{\sqrt 2}g(\frac x 2)\widehat{\phi}(\frac x 2).$$ is an
orthonormal wavelet.
\end{theorem}

  A good introductory
description of these constructions appears in \cite{Str}.

The theorem's requirements that $h$ satisfy $|h(0)|=\sqrt2$ and that $h$
be in $C^1$ are natural restrictions in order to make
the infinite product converge.  The condition that $h$ is
nonvanishing on $[-\frac 14,\frac 14]$ appears less natural; it is
used in the proof to ensure $L^2$ convergence of the infinite
product and thus the orthonormality of the translates of $\phi$. A
famous example due to A. Cohen \cite{Coh} showed that this 
nonvanishing condition cannot be entirely removed, as its removal can
lead to functions $\phi$ and $\psi$ whose translates are not
orthonormal. Cohen took $h=\frac{1+e^{-6\pi i x}}{\sqrt 2}$, which
resulted in a stretched out version of the Haar scaling function and
wavelet, $\phi=\frac 13 \chi_{[0,3)}$ and $\psi=\frac 13(\chi_{[0,\frac 32)}-\chi_{[\frac 32, 3)}).$
Cohen's $\psi$ is not a wavelet, since its translates are not orthonormal.
However, the classical Theorem \ref{classical} can be extended to
accommodate this and similar examples if we allow (normalized tight)
frame wavelets.
The following generalization of Theorem \ref{classical}
was first proven by Lawton in
\cite{Law}, with a generalized form appearing later in Bratteli and Jorgensen's book \cite{BJ}:

\begin{theorem} \label{Palle}
Suppose $h,\,g$ are periodic Lipschitz continuous function in $L^2(\mathbb R)$, which satisfy
$|h(0)|=\sqrt 2$ and the filter equations (\ref{orth1}), (\ref{orth2}) and (\ref{orth3}).
Then the construction $\widehat{\phi}(x)=\prod_{j=1}^{\infty} \frac 1{\sqrt 2}h(2^{-j}x)$
 produces an $L^2$ function $\phi$ (whose translates are not necessarily
orthogonal), and the function $\widehat{\psi}(x)=\frac
1{\sqrt 2}g(\frac x 2)\widehat{\phi}(\frac x 2)$ is a (normalized tight) frame
wavelet for dilation by $2$ in $L^2(\mathbb{R})$.
\end{theorem}

The classical filter construction that was used in both the original theorem
and this extension was generalized to GMRA's in
\cite{JKL}.  Just as in the classical case, it is natural to
begin by building filters from wavelets, and then determine
conditions under which the process can be reversed. Accordingly,
suppose first that $\psi$ is a non-MRA orthonormal wavelet for
dilation by $2$ in $L^2(\mathbb R)$ whose multiplicity function is
bounded by $c$.   We have seen in Section 2 that
there then exist generalized
scaling functions $\phi_1,\cdots,\phi_c$. Since $\widehat{\phi}$ and $\widehat {\psi}$ are in $\widehat{V}_{1}$,
we can show (see \cite{JKL}) there exist
periodic functions $h_{i,j}$ and $g_{j}$, supported on the
periodization of $S_j$, such that
\begin{equation}
\widehat{\phi_i}(x)=\frac 1{\sqrt{2}}\sum_{j=1}^c h_{i,j}(\frac x2)\widehat{\phi_j}(\frac x2)\label{gen1}
\end{equation}
and
\begin{equation}
\widehat{\psi}(x)=\frac 1{\sqrt{2}}\sum_{j=1}^c g_{j}(\frac x2)\widehat{\phi_j}(\frac x2).\label{gen2}
\end{equation}
These \emph {generalized filters} $g_{j}$ and $h_{i,j}$ satisfy
orthonormality-like conditions that are generalizations of the
classical conditions (\ref{orth1}), (\ref{orth2}) and
(\ref{orth3}):
\begin{equation}
\sum_{j=1}^c h_{i,j}(\frac x2)\overline {h_{k,j}(\frac x2)}+h_{i,j}(\frac {x+1}2)\overline {h_{k,j}(\frac {x+1}2)}=2\delta_{i,k}\chi_{S_i}(x),\label{gen3}
\end{equation}
\begin{equation}
\sum_{j=1}^c |g_{j}(\frac x2)|^2+|g_{j}(\frac {x+1}2)|^2=2,\label{gen4}
\end{equation}
 and
\begin{equation}
\sum_{j=1}^c h_{i,j}(\frac x2)\overline {g_{j}(\frac x2)}+ h_{i,j}(\frac {x+1}2)\overline {g_{j}(\frac {x+1}2)}=0,\label{gen5}
\end{equation}

For example, for the Journ\'e wavelet these filters are defined (see \cite{Cou}) by:

$$h_{1,1}=\sqrt 2\chi_{[-\frac 27,-\frac 14)\cup (-\frac 17,\frac
17)\cup [\frac 14,\frac 27)}$$

$$h_{1,2}=0,$$

$$h_{2,1}=\sqrt 2\chi_{[-\frac 47,-\frac 12)\cup [\frac 12,\frac 47)}$$

$$h_{2,2}=0$$

$$g_1=\sqrt 2\chi_{[-\frac 14,-\frac 17)\cup [\frac 17,\frac 14)}$$

$$g_2=\sqrt 2\chi_{[-\frac 17,\frac 17)}.$$

To use generalized filters to build new wavelets, we
reverse this procedure, just as in the classical case. In
order to first build filters, we use functions on the disjoint
union of the $S_j$'s whose values are $\sqrt{2}$ times
unitary matrices, with different dimensions for different values
of $x$. We need the values of the filters to be $\sqrt{2}$
times unitary matrices in order to satisfy the generalized
orthonormality conditions (\ref{gen3}), (\ref{gen4}), and
(\ref{gen5}). The matrices of filter values have different
dimensions depending on how many of the sets $S_j$ the point $x$
and its preimages are in. Once we have the filters, we build the
generalized scaling function using an infinite product of matrices
that comes from the iteration of Eq.(\ref{gen1}). The wavelet
is then produced by Eq.(\ref{gen2}). Conditions that make
this possible are described in the following generalization (see \cite{BJMP})
of the Lawton theorem :

\begin{theorem} \label{bjmp}
Given a multiplicity function $m$ for a GMRA, suppose $\{h_{i,j}\}$ and
$\{g_{j}\}$ are periodic functions that are supported on the
periodization of $S_j=\{x: m(x)\geq j\}$, Lipschitz continuous in a neighborhood of
the origin, and that satisfy the three generalized orthonormality conditions
(\ref{gen3}),(\ref{gen4}), and (\ref{gen5}).  Suppose in addition the $h_{i,j}$ satisfy the
generalized lowpass condition
$|h_{i,j}(0)|=\sqrt{2}\delta_{(i,1)}\delta_{(j,1)}$.
Write $H$ for the matrix $(h_{i,j})$.  Then the components of\newline
$\prod_{k=1}^{\infty} \frac 1{\sqrt 2}H(2^{-k}x)$
converge pointwise to $L^2$ functions.  If we let
$\{\widehat{\phi_i}\}_{i=1}^c$ be the first column of this product, then the translates of $\{\phi_i\}$ determine
the core subspace $V_0$ of a GMRA, and
 $$\widehat{\psi}(x )
\equiv\frac 1{\sqrt 2}g_{1}(\frac x2 )\widehat{ \phi_1}(\frac x2 )+g_{2}(\frac x2 )\widehat{ \phi_2}(\frac x2 )$$
is the Fourier transform of a normalized tight frame wavelet
on $L^2(\mathbb R).$
\end{theorem}

The proof, like that of \cite{BJ}, proceeds by using matrices of
values of the filters to define partial isometries that satisfy
relations similar to those defining a  Cuntz algebra as in \cite{Cun}.
We will use this theorem in the next section to build our main example.

\section{The Example}
For any fixed positive integer $r$, we will now use Theorem \ref{bjmp} to build  an example of a $C^r$
frame wavelet with
$C^{\infty}$ Fourier transform, whose multiplicity function is that
of Journ\'e:
$$m(x) = \left\{\begin{array}{ll}
2 & x\in[-\frac 17,\frac 17) \\
1 & x\in\pm[\frac 17,\frac 27)\cup \pm [\frac 37,\frac 12)\\
  0 &
\mbox{otherwise}
\end{array}\right.$$
To do this, we will first choose low-pass filters that force the functions $\widehat{\phi_1}$ and $\widehat{\phi_2}$
in the infinite matrix product to both be smooth and have rapid decay.  We will then define the high-pass filters, and show that
the resulting wavelet inherits the properties we seek.

{}From the equation for $m$, we see that the sets $S_j=\{x:m(x)\geq j\}$ used in the filter construction are
here given by $S_1=[-\frac 27,\frac 27)\cup\pm[\frac 37,\frac 12)$ and $S_2=[-\frac 17, \frac 17)$.
We will define the filter functions $h_{i,j}$, $1\leq i,j\leq 2$
in terms of a
classical filter $p$.  Accordingly, with $r$ fixed as above, let $p$ be a periodic (with period 1), real-valued
function that satisfies:
\begin{enumerate}
\item $p\in C^{\infty}(\mathbb R)$
\item $p(x)=p(-x)$
 \item $|p(x)|^2+|p(x+\frac12)|^2=2$
 \item $p^{(k)}(x)=0$ for $x\in[\frac17,\frac3{14}]\cup[\frac37,\frac12]$, $k\geq 0$.
\item $|p^{(r+2)}(x)|<1$ when $|x-\frac 3{14}|<\frac 1{112}$ or $|x-\frac 17|<\frac 1{112}$
\end{enumerate}
To see that such a function exists, note that it suffices to
define $p\in C^{\infty}([\frac3{28},\frac14]\cup[\frac{11}{28},\frac12])$, satisfying conditions (4) and (5), and with
$p(\frac3{28})=p(\frac14)=p(\frac{11}{28})=1$ and $p^{(k)}(\frac3{28})=p^{(k)}(\frac14)=p^{(k)}(\frac{11}{28})=0$ for $k\geq 1$.
We can then define $p$ on $[0,\frac3{28}]\cup[\frac14,\frac{11}{28}]$ by
$p(x)=\sqrt{2-p(\frac12-x)^2}$.  We then extend $p$ to $[-\frac12,0]$ using
the requirement that $p$ be symmetric, and finally extend $p$ to the whole line using periodicity.
Figure 1 below shows a $C^1$ approximation to $p.$

\begin{figure}[ht]
\setlength{\unitlength}{300bp}
\begin{picture}(1.04,0.6489315)(-0.04,-0.032446575)
\put(0,0){\includegraphics[width=\unitlength]{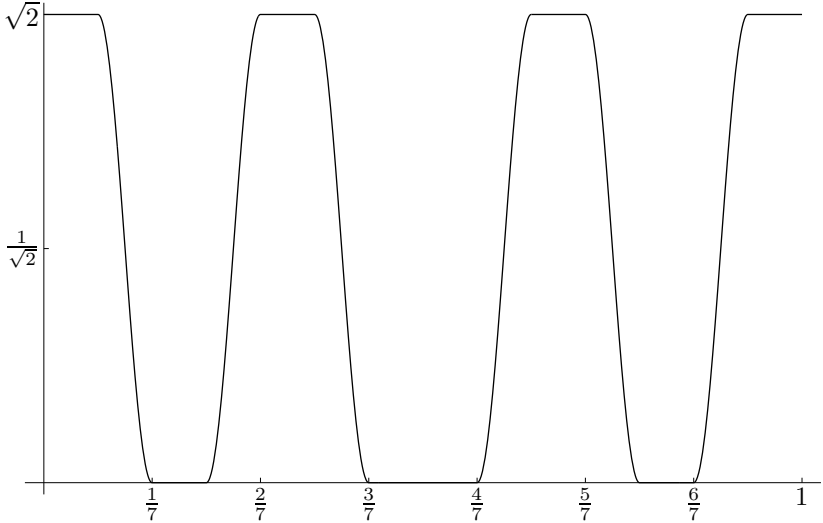}}
\put(0.15986,0.00889){\makebox(0,0)[t]{$\frac{1}{7}$}}
\put(0.29592,0.00889){\makebox(0,0)[t]{$\frac{2}{7}$}}
\put(0.43197,0.00889){\makebox(0,0)[t]{$\frac{3}{7}$}}
\put(0.56803,0.00889){\makebox(0,0)[t]{$\frac{4}{7}$}}
\put(0.70408,0.00889){\makebox(0,0)[t]{$\frac{5}{7}$}}
\put(0.84014,0.00889){\makebox(0,0)[t]{$\frac{6}{7}$}}
\put(0.97619,0.00889){\makebox(0,0)[t]{$1$}}
\put(0.01798,0.30902){\makebox(0,0)[r]{$\frac{1}{\sqrt{2}}$}}
\put(0.01798,0.60332){\makebox(0,0)[r]{$\sqrt{2}$}}
\end{picture}
\caption{A $C^1$ approximation to $p$}
\label{Fig1}
\end{figure}

 We now use this function $p$ to build the low-pass filters $\{h_{i,j}\}$:
\begin{equation}\label{h11} h_{1,1}(x)=\begin{cases} p(x)&x\in [-\frac 2{7},\frac 27)\cr
0&\text{\rm \ otherwise}\end{cases}
\end{equation}

\begin{equation}\label{h12}h_{1,2}(x)=\begin{cases} p(x+\frac 12)&x\in [-\frac 1{7},\frac 17)\cr
0&\text{\rm \ otherwise}\end{cases}
\end{equation}

\begin{equation}\label{h21}h_{2,1}(x)=\begin{cases} \sqrt 2 &x\in \pm [\frac 37,\frac 12)\cr
0&\text{\rm \ otherwise}\end{cases}
\end{equation}

\begin{equation}\label{h22}h_{2,2}(x)=0
\end{equation}

Figures 2 through 4 show $C^1$ approximations to the three nonzero low-pass filters.

\begin{figure}[ht]
\setlength{\unitlength}{300bp}
\begin{picture}(1.04,0.6489315)(-0.04,-0.032446575)
\put(0,0){\includegraphics[width=\unitlength]{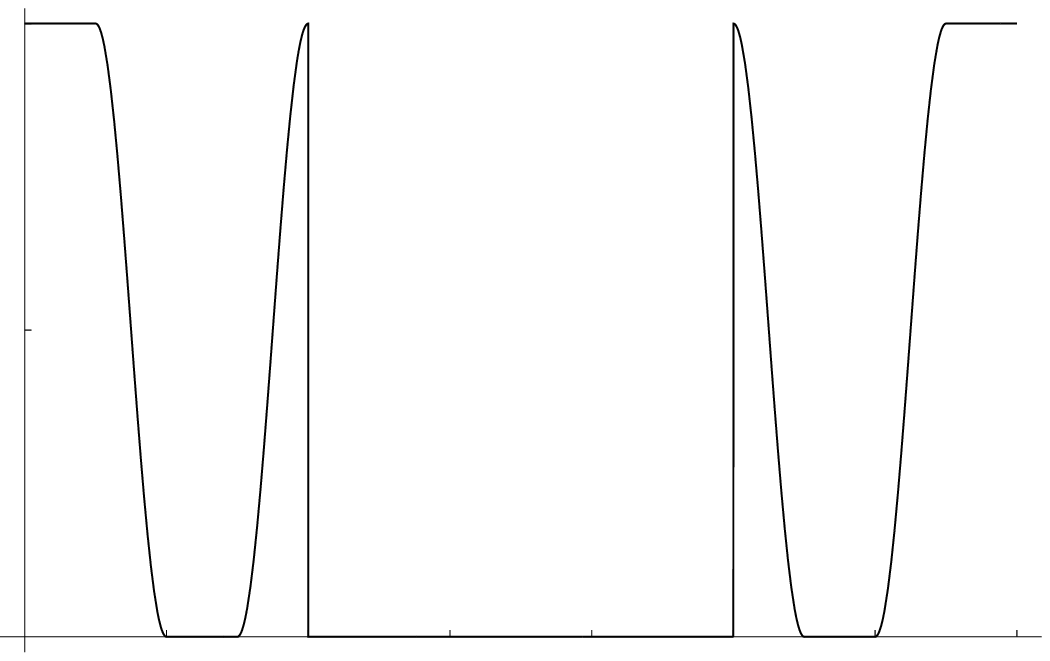}}
\put(0.15986,0.00889){\makebox(0,0)[t]{$\frac{1}{7}$}}
\put(0.29592,0.00889){\makebox(0,0)[t]{$\frac{2}{7}$}}
\put(0.43197,0.00889){\makebox(0,0)[t]{$\frac{3}{7}$}}
\put(0.56803,0.00889){\makebox(0,0)[t]{$\frac{4}{7}$}}
\put(0.70408,0.00889){\makebox(0,0)[t]{$\frac{5}{7}$}}
\put(0.84014,0.00889){\makebox(0,0)[t]{$\frac{6}{7}$}}
\put(0.97619,0.00889){\makebox(0,0)[t]{$1$}}
\put(0.01798,0.30902){\makebox(0,0)[r]{$\frac{1}{\sqrt{2}}$}}
\put(0.01798,0.60332){\makebox(0,0)[r]{$\sqrt{2}$}}
\end{picture}
\caption{An approximation to $h_{1,1}$}
\label{Fig2}
\end{figure}

\begin{figure}[ht]
\setlength{\unitlength}{300bp}
\begin{picture}(1.04,0.6489315)(-0.04,-0.032446575)
\put(0,0){\includegraphics[width=\unitlength]{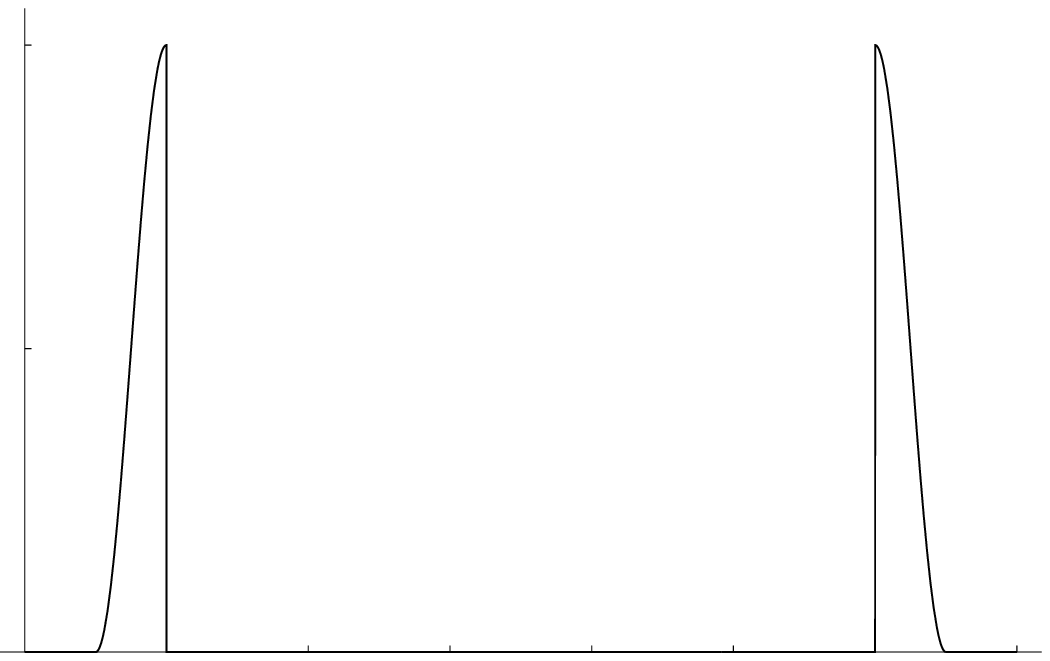}}
\put(0.15986,0){\makebox(0,0)[t]{$\frac{1}{7}$}}
\put(0.29592,0){\makebox(0,0)[t]{$\frac{2}{7}$}}
\put(0.43197,0){\makebox(0,0)[t]{$\frac{3}{7}$}}
\put(0.56803,0){\makebox(0,0)[t]{$\frac{4}{7}$}}
\put(0.70408,0){\makebox(0,0)[t]{$\frac{5}{7}$}}
\put(0.84014,0){\makebox(0,0)[t]{$\frac{6}{7}$}}
\put(0.97619,0){\makebox(0,0)[t]{$1$}}
\put(0.01798,0.3){\makebox(0,0)[r]{$\frac{1}{\sqrt{2}}$}}
\put(0.01798,0.59){\makebox(0,0)[r]{$\sqrt{2}$}}
\end{picture}
\caption{An approximation to $h_{1,2}$}
\label{Fig3}
\end{figure}

\begin{figure}[ht]
\setlength{\unitlength}{300bp}
\begin{picture}(1.04,0.6489315)(-0.04,-0.032446575)
\put(0,0){\includegraphics[width=\unitlength]{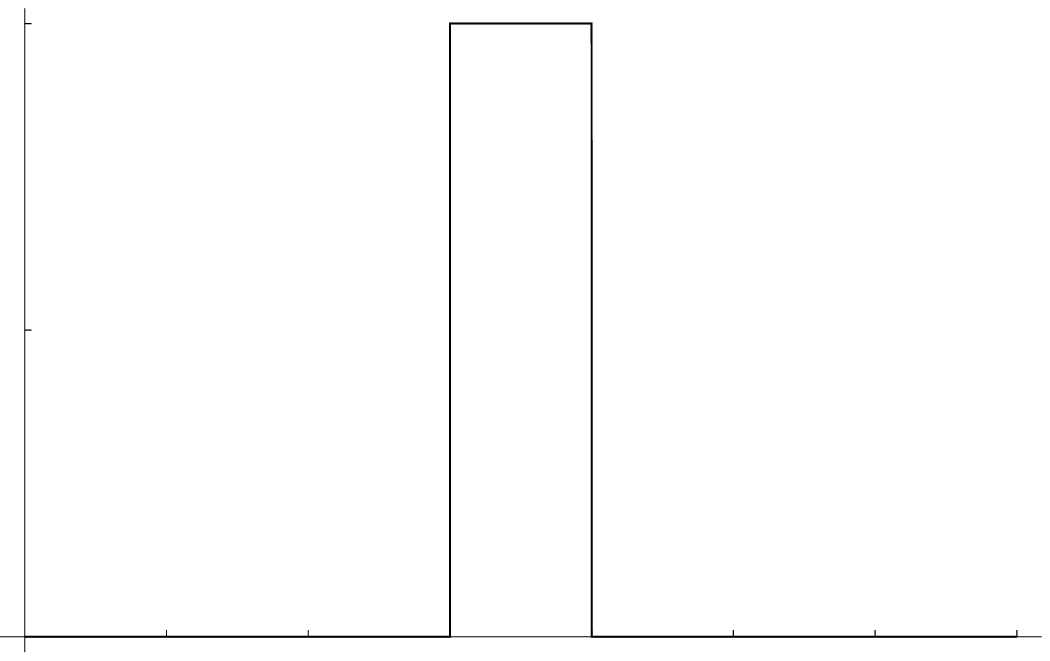}}
\put(0.15986,0.0088){\makebox(0,0)[t]{$\frac{1}{7}$}}
\put(0.29592,0.0088){\makebox(0,0)[t]{$\frac{2}{7}$}}
\put(0.43197,0.0088){\makebox(0,0)[t]{$\frac{3}{7}$}}
\put(0.56803,0.0088){\makebox(0,0)[t]{$\frac{4}{7}$}}
\put(0.70408,0.0088){\makebox(0,0)[t]{$\frac{5}{7}$}}
\put(0.84014,0.0088){\makebox(0,0)[t]{$\frac{6}{7}$}}
\put(0.97619,0.0088){\makebox(0,0)[t]{$1$}}
\put(0.01798,0.30902){\makebox(0,0)[r]{$\frac{1}{\sqrt{2}}$}}
\put(0.01798,0.60332){\makebox(0,0)[r]{$\sqrt{2}$}}
\end{picture}
\caption{An approximation to $h_{2,1}$}
\label{Fig4}
\end{figure}

It is easy to check that these functions satisfy the filter equation (\ref{gen3}).
For example, to check that
$$|h_{1,1}(\frac x2)|^2+|h_{1,2}(\frac x2)|^2+|h_{1,1}(\frac{x+1}2)|^2+|h_{1,2}(\frac{x+1}2)|^2=2\chi_{[\frac{-2}7,\frac27]\cup\pm[\frac37, \frac 12]}(x),$$
we note that the left hand side gives
$$\left\{\begin{array}{lcl}|p(\frac x2)|^2+|p(\frac{x+1}2)|^2+0+0&=&2\mbox{ for } x\in[\frac{-2}7,\frac27]\\
|p(\frac x2)|^2+0+0+0&=&0\mbox{ for }x\in\pm[\frac27,\frac 37]\\
|p(\frac x2)|^2+0+|p(\frac{x+1}2)|^2+0&=&2\mbox{ for }x\in\pm[\frac37,\frac12]\end{array}\right.$$
  The other cases of (\ref{gen3}) follow
similarly from the definitions of the $h_{i,j}$ and the properties of $p$.  It is also easy to check that the $\{h_{i,j}\}$ have proper support and satisfy the low-pass condition, $h_{1,1}(0)=\sqrt{2}$ and $h_{i,j}(0)=0$ for $(i,j)\neq(0,0).$

As in Theorem \ref{bjmp}, we will define the functions $\widehat{\phi_1}$ and $\widehat{\phi_2}$ (and thus the GMRA that they determine) in terms of the low-pass filter $\{h_{i,j}\}$.  First we will need the following technical lemma:

\begin{lemma} \label{partialproduct}
For every integer $n\geq 1$, the $(3n)^{th}$ partial product of the filter matrices $\prod_{j=1}^{3n}\frac 1{\sqrt 2} \left(\begin{matrix} h_{1,1}(\frac x{2^j})&h_{1,2}(\frac x{2^j})\\
h_{2,1}(\frac x{2^j})&0\end{matrix} \right)$ is a lower triangular matrix $\left(\begin{matrix} a_n(x)&0\\
c_{n}(x)&d_n(x)\end{matrix} \right)$ whose nonzero components have the following support:

$$\text{support}(a_n)=\left(8^n\mathbb{Z}+\left([-\frac27,\frac27]\cup\pm[\frac37,\frac47]
\bigcup\cup_{k=0}^{n-2}A_{k}\bigcup\cup_{k=0}^{n-1}B_k\right)\right)$$

$$\bigcup\left(4\cdot8^{n-1}\mathbb{Z}+A_{n-1}\right)$$

$$\text{support}(c_n)=8^n\mathbb{Z}+\cup_{k=0}^{n-1}C_k$$

$$\text{support}(d_n)=8^n\mathbb{Z}+C_{n-1}$$

\noindent where
$$A_k=\bigcup_{\{a_j\}\in\{1,-1\}^{k+1}}\sum_{j=0}^k 2 a_j 8^j+\pm[\frac17,\frac27],$$

$$B_k=\bigcup_{\{a_j\}\in\{1,-1\}^{k+1}}\sum_{j=0}^k 4 a_j 8^j+\pm[\frac37,\frac47],$$

$$\text{and}\quad C_k=\bigcup_{\{a_j\}\in\{1,-1\}^{k+1}}\sum_{j=0}^k a_j 8^j+ [-\frac17,\frac17].$$

\end{lemma}

\begin{proof}
The proof is by induction on $n.$  Both the $n=1$ case and the induction step follow directly from
the observation that
$$\prod_{j=1}^{3}\left(\begin{matrix} h_{1,1}(\frac x{2^j})&h_{1,2}(\frac x{2^j})\\
h_{2,1}(\frac x{2^j})&0\end{matrix} \right)=\left(\begin{matrix} t_{1,1}(x)&0\\
t_{2,1}(x)&t_{2,2}(x)\end{matrix} \right)$$
where $$t_{1,1}(x)=h_{1,1}(\frac x{2})h_{1,1}(\frac x{4})h_{1,1}(\frac x{8})+h_{1,2}(\frac x{2})h_{2,1}(\frac x{4})h_{1,1}(\frac x{8})
+h_{1,1}(\frac x{2})h_{1,2}(\frac x{4})h_{2,1}(\frac x{8})$$ has support $$8\mathbb{Z}+([-\frac27,\frac27]\cup\pm[\frac37,\frac47])
\bigcup 4\mathbb{Z}+2\pm[\frac17,\frac27]\bigcup 8\mathbb{Z}+4\pm[\frac37,\frac47],$$ and both
$$t_{2,1}(x)=h_{2,1}(\frac x{2})h_{1,1}(\frac x{4})h_{1,1}(\frac x{8})$$ and $$t_{2,2}(x)=h_{2,1}(\frac x{2})h_{1,1}(\frac x{4})h_{1,2}(\frac x{8})$$
have support $$8\mathbb{Z}\pm[\frac67,\frac87].$$
\end{proof}

\begin{remark}
Since \cite{BJMP} shows that the infinite product of filter matrices, $$\prod_{j=1}^\infty\frac 1{\sqrt 2} \left(\begin{matrix} h_{1,1}(\frac x{2^j})&h_{1,2}(\frac x{2^j})\\
h_{2,1}(\frac x{2^j})&0\end{matrix} \right)$$ converges pointwise, Lemma \ref{partialproduct} gives the support of the elements of this matrix.  In
particular, the support of the upper left element (which will be used to define $\widehat{\phi_1}$) is $[-\frac27,\frac27]\cup\pm[\frac37,\frac47]\bigcup\cup_{k=0}^{\infty}A_{k}\bigcup\cup_{k=0}^{\infty}B_k,$
and the support of the lower left element (which will be used to define $\widehat{\phi_2}$) is $\cup_{k=0}^{\infty}C_k$.  Note that these supports are disjoint.
\end{remark}

We will use this result to characterize the GMRA for the wavelet we will build.

\begin{proposition}
\label{GMRAex}
The functions $\widehat{\phi_1}$ and $\widehat{\phi_2}$ given by the upper left and lower left elements of the infinite product
matrix $\prod_{j=1}^\infty\frac 1{\sqrt 2} \left(\begin{matrix} h_{1,1}(\frac x{2^j})&h_{1,2}(\frac x{2^j})\\
h_{2,1}(\frac x{2^j})&0\end{matrix} \right)$ determine a GMRA with multiplicity function the same as the Journ\'e
wavelet.
\end{proposition}

\begin{proof}
Define $\widehat{V_0}$ to be the closed linear span of exponentials times $\widehat{\phi_1}$ and $\widehat{\phi_2}$, or in other words, define  $V_0$ to be  the closed linear span of translates of $\phi_1$ and $\phi_2$.  A GMRA can be built from
this core subspace $V_0$ in a natural way, by letting $V_j=\{\delta^j(f)\equiv\sqrt{2}^jf(2^jx)\}_{f\in
V_0}$.  By \cite{BJMP}, the resulting structure does satisfy all the components of Definition \ref{GMRA}.   To find the multiplicity function associated with this GMRA, we note first that Lemma \ref{partialproduct} shows that the supports of
$\widehat{\phi_1}$ and $\widehat{\phi_2}$ are disjoint, so that the multiplicity function is the sum of the multiplicity functions determined by the two functions separately.  By \cite{Ba}, the multiplicity function associated with
$\widehat{\phi_j}$ is the characteristic function of the support of the periodic function $\text{Per}\,\phi_j(x)=\sum_{l\in\mathbb{Z}}|\widehat{\phi_j}(x+l)|^2$.  From Lemma \ref{partialproduct}, we see that these sets are indeed the Journe sets $S_1$ for $\widehat{\phi_1}$ and $S_2$ for $\widehat{\phi_2}$.
\end{proof}

\begin{remark} Neither $\text{Per}\,\phi_1$ nor $\text{Per}\,\phi_2$ is  bounded away from 0 on its support, so
that by \cite{BLi}, the translates of $\phi_1$ and $\phi_2$ do not even form a frame for
their closed linear span $V_0$.  (This is true in the Cohen example as well.)
\end{remark}
Figures \ref{Fig5} and \ref{Fig6} below show $C^1$ approximations to the graphs of $\widehat{\phi_1}$ and $\widehat{\phi_2}$
near the origin.  As shown in Lemma \ref{partialproduct}, both functions have unbounded support, but the portions of these functions that do not appear in Figure \ref{Fig5} and Figure \ref{Fig6} have values less than .01 and .002 respectively.  We will use $\widehat{\phi}_1$ and $\widehat{\phi}_2$ to define our $C^r$ frame wavelet with $C^{\infty}$ Fourier transform.
To do this we will need the following properties of $\widehat{\phi}_1$ and $\widehat{\phi}_2$:

\begin{figure}[ht]
\setlength{\unitlength}{300bp}
\begin{picture}(1.04,0.6489315)(-0.04,-0.032446575)
\put(0,0){\includegraphics[width=\unitlength]{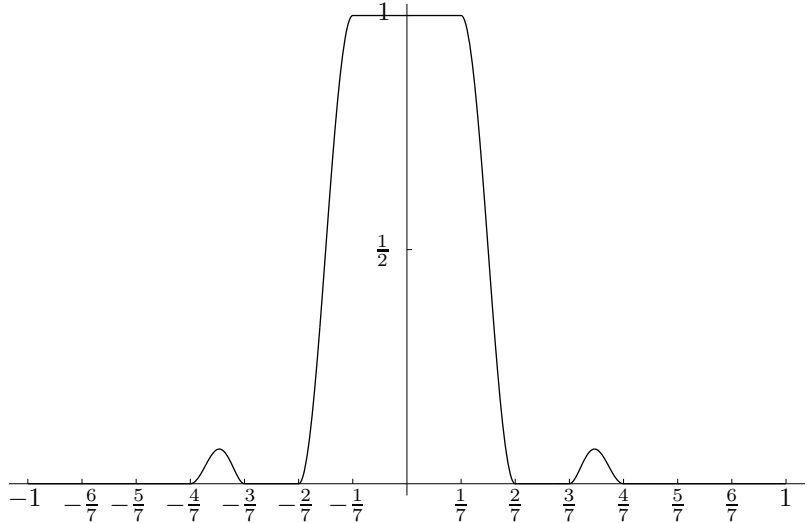}}
\put(0.02,0.009){\makebox(0,0)[t]{$-1$}}
\put(0.091,0.00889){\makebox(0,0)[t]{$-\frac{6}{7}$}}
\put(0.15,0.00889){\makebox(0,0)[t]{$-\frac{5}{7}$}}
\put(0.22,0.00889){\makebox(0,0)[t]{$-\frac{4}{7}$}}
\put(0.29,0.00889){\makebox(0,0)[t]{$-\frac{3}{7}$}}
\put(0.36,0.00889){\makebox(0,0)[t]{$-\frac{2}{7}$}}
\put(0.425,0.00889){\makebox(0,0)[t]{$-\frac{1}{7}$}}
\put(0.56803,0.00889){\makebox(0,0)[t]{$\frac{1}{7}$}}
\put(0.63606,0.00889){\makebox(0,0)[t]{$\frac{2}{7}$}}
\put(0.70408,0.00889){\makebox(0,0)[t]{$\frac{3}{7}$}}
\put(0.77211,0.00889){\makebox(0,0)[t]{$\frac{4}{7}$}}
\put(0.84014,0.00889){\makebox(0,0)[t]{$\frac{5}{7}$}}
\put(0.90817,0.00889){\makebox(0,0)[t]{$\frac{6}{7}$}}
\put(0.97619,0.00889){\makebox(0,0)[t]{$1$}}
\put(0.48,0.30902){\makebox(0,0)[r]{$\frac{1}{2}$}}
\put(0.48,0.61){\makebox(0,0)[r]{$1$}}

\end{picture}
\caption{A $C^1$ approximation to $\widehat{\phi_1}$}
\label{Fig5}
\end{figure}

\begin{figure}[ht]
\setlength{\unitlength}{300bp}
\begin{picture}(1.04,0.6489315)(-0.04,-0.032446575)
\put(0,0){\includegraphics[width=\unitlength]{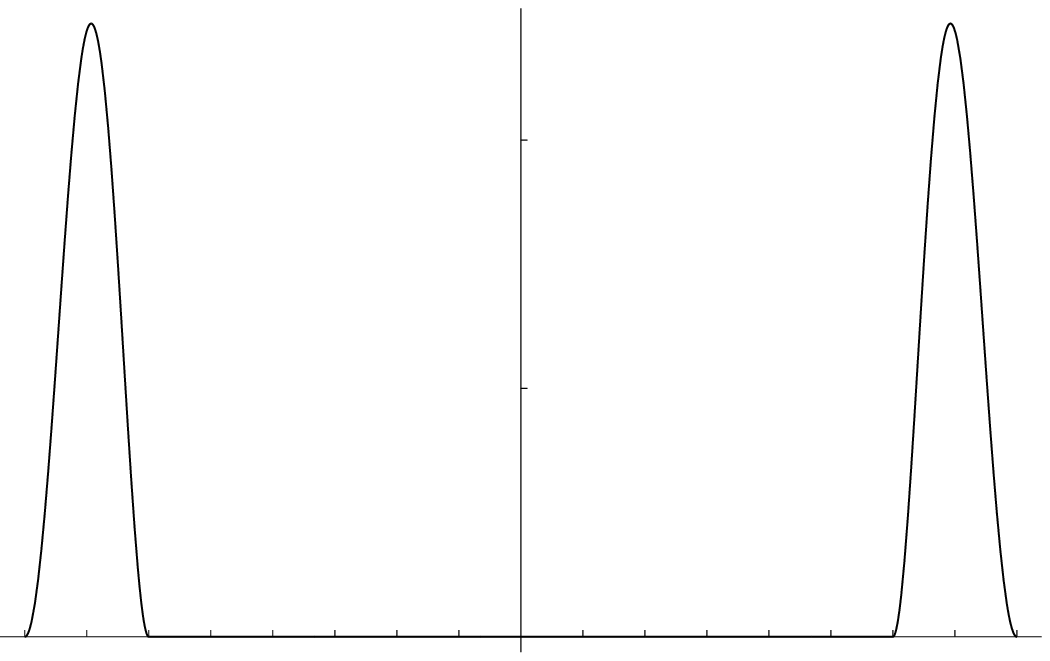}}
\put(0.02,0.009){\makebox(0,0)[t]{$-\frac87$}}
\put(0.14,0.00889){\makebox(0,0)[t]{$-\frac{6}{7}$}}
\put(0.26,0.00889){\makebox(0,0)[t]{$-\frac{4}{7}$}}
\put(0.38,0.00889){\makebox(0,0)[t]{$-\frac{2}{7}$}}
\put(0.62,0.00889){\makebox(0,0)[t]{$\frac{2}{7}$}}
\put(0.74,0.00889){\makebox(0,0)[t]{$\frac{4}{7}$}}
\put(0.86,0.00889){\makebox(0,0)[t]{$\frac{6}{7}$}}
\put(0.98,0.00889){\makebox(0,0)[t]{$\frac87$}}
\put(0.48,0.25){\makebox(0,0)[r]{$.03$}}
\put(0.48,0.5){\makebox(0,0)[r]{$.06$}}

\end{picture}
\caption{A $C^1$ approximation to $\widehat{\phi_2}$}
\label{Fig6}
\end{figure}

\begin{proposition}\label{Cinfty}The functions $\widehat{\phi}_1$ and $\widehat{\phi}_2$ defined in Proposition \ref{GMRAex} are
$C^{\infty}$.
\end{proposition}

\begin{proof}We fix an $a$ and show that $\widehat{\phi}_1$ and $\widehat{\phi}_2$
are $C^{\infty}$ in a neighborhood of $a$.
Because $h_{2,1}$ is 0 on $(-\frac 37, \frac 37)$, all but a finite
number of the lower triangular products of three successive factors described in Lemma \ref{partialproduct} are
in fact diagonal with single terms on the diagonal.  Thus, each of the entries in the first column of the infinite product matrix is the sum
 of only a finite number of terms.  Since $h_{1,1}=\sqrt{2}$ in a neighborhood of $0$, each of the terms has only a finite number of factors not equal to 1.
It will suffice to show that each of
these finite products is $C^{\infty}$.

By construction, we have that the $h_{i,j}$ are all
$C^\infty$ everywhere except for $h_{1,1}$ at $n\pm\frac 27$, $h_{1,2}$ at $n\pm\frac 17$, and
$h_{2,1}$ at $n\pm\frac 37$ for $n\in\Bbb Z$. We will show that whenever one of these discontinuities occurs as a factor
in the infinite product, it is cancelled by a following factor that is $0$ at the point of discontinuity.
{}From the formula for the $t_{i,j}$ given in Lemma \ref{partialproduct}, we see that any term in the infinite product that contains a factor of $h_{2,1}(n\pm\frac 37)$ must also contain
a factor of one of the forms $h_{1,1}(n\pm\frac 2{7})$, $h_{1,1}(n\pm\frac 3{14})$, $h_{1,2}(n\pm\frac 3{14})$, or $h_{1,2}(n\pm\frac 2{7})$.
The last three possibilities are $0$ in a neighborhood of the point in question, so if we have a discontinuous factor of $h_{2,1}$, it is
either cancelled out by a $0$ factor, or we also have a factor of $h_{1,1}(n\pm\frac 27)$ with a smaller $n$.  Similarly,
any term in the infinite product that contains a factor of $h_{1,1}(n\pm\frac 27)$ must also contain
a factor of one of the forms $h_{1,1}(n\pm\frac 5{14})$,  $h_{1,1}(n\pm\frac 1{7})$,
$h_{1,2}(n\pm\frac 5{14})$, or $h_{1,2}(n\pm\frac 1{7})$. The first three of these possibilities are 0 at the points in question,
 so any discontinuous factor of $h_{1,1}$ is either cancelled out by a $0$ factor, or is followed by a factor of $h_{1,2}(n\pm\frac 17)$ with an equal or smaller
$n$.
Finally, any term in the infinite product that contains a factor
of $h_{1,2}(n\pm\frac 17)$ must also contain
a factor of either the form $h_{2,1}(n\pm\frac 1{14})$ or the form $h_{2,1}(n\pm\frac 3{7})$. The first of these possibilities is $0$
at the point in question; the second possibility throws us back into the first type of discontinuity we considered above, but with a
smaller $n$.  We can repeat the above sequence of arguments if necessary, noting that each succeeding factor is evaluated at a point half the distance from
the origin as its predecessor, so that the chain above must eventually end with a factor of $h_{1,1}$ and thus with $0$.

Thus each of the finite products that occur as terms in $\widehat{\phi}_1$ and $\widehat{\phi}_2$ has the form of product $f(x)=f_1(x)f_2(x)$ in
a neighborhood of $a$, where $f_1(x)$ is $C^{\infty}$ on $(a-\epsilon,a+\epsilon)$, with its value and all its derivatives, $f_1^{(n)}(a)=0$, $n\geq 0$,
and $f_2(x)$ is $C^{\infty}$ on $(a-\epsilon,a)\cup(a,a+\epsilon)$, with a jump discontinuity at $a$, but with all derivatives satisfying
$\lim_{x\rightarrow a^-}f_2^{(n)}(x)=\lim_{x\rightarrow a^+}f_2^{(n)}(x)=0$, $n\geq 1$.  Thus, $f$ is $C^{\infty}$ on $(a-\epsilon,a+\epsilon)$ with
$f^{(n)}(a)=0$, $n\geq 0$.
\end{proof}

The next result shows that the $\widehat{\phi_i}$ vanish rapidly at $\infty,$ to an extent determined by the integer $r$ 
chosen at the beginning of the example.

\begin{proposition}\label{xr}The functions $\widehat{\phi}_1$ and $\widehat{\phi}_2$ defined in Proposition \ref{GMRAex}  satisfy
 $|x|^{r+1}\widehat{\phi_i}(x)\rightarrow 0$ as $x\rightarrow\infty.$
\end{proposition}

\begin{proof}
Using Lemma \ref{partialproduct} we see that it will suffice to show that given an $\epsilon>0$, there
exists $N$ such that $|x|^{r+1}|\widehat{\phi_i}(x)|<\epsilon$ for $x\in A_n\cup B_n\cup C_n$ with $n>N$.
Recall from that lemma that $x\in A_n\cup B_n\cup C_n$ has the form $x=d\left(\sum_{j=0}^{n} a_{j}8^j+\pm r_0\right)$,
where $d\in{1,2,4}$, $a_j=\pm1$, and where $0<r_0<\frac17.$
For any such $x$,  $\widehat{\phi}_i(x)$ is the sum of
less than $3^{n+1}$ terms of the form $\prod_{l=1}^{\infty}\frac 1{\sqrt{2}}h_{i_l,j_l}(\frac x{2^l})$, since once $x<\frac 67,$ we
have $h_{2,1}(x)=0$, so that the triangular matrix described in Lemma \ref{partialproduct} is actually diagonal with single term elements.
To see how the argument of the $h_{i_l,j_l}$ changes as we move from $\frac x{2^k}$ to $\frac x{2^{k+3}}$, we note that if $y$ of the form $y=\sum_{j=0}^{n-k} a_{k+j}8^j+a_{k-1}r_k$ with
$\frac 17-\delta<r_0<\frac17$,  then $\frac y2=4\left(\sum_{j=0}^{n-k-1}a_{k+j+1}8^j+a_k r_{k+1}\right)$, $\frac y4=2\left(\sum_{j=0}^{n-k-1}a_{k+j+1}8^j+a_k r_{k+1}\right)$, and $\frac y8=\sum_{j=0}^{n-k-1}a_{k+j+1}8^j+a_k r_{k+1}$.  The remainder terms, $d\cdot r_{k+1}$, take one of two forms, depending on the value of $a_{k-1}$. If $a_{k-1}=+1$, then $4r_{k+1}$ is within $\frac{\delta}2$ of $\frac47$, so that $\frac17-\frac{\delta} 8<r_{k+1}<\frac17$, while if $a_{k-1}=-1$, $4r_{k+1}$ is within $\frac{\delta}2$ of $\frac37$, so that $\frac3{14}<2r_{k+1}<\frac3{14}+\frac{\delta}4$, and
$\frac17-\frac1{28}<r_{k+1}<\frac17$. In the first case we will make use of the fact that as $y$ decreases by a factor of $\frac18$, $\delta$ decreases by a factor of $\frac18$.  In the second case we make use of the factor $h_{i,j}(z)$ with $\|z-\frac 3{14}\|<\frac{\delta}4$, where $\|\cdot\|$ denotes the distance from the nearest integer.  We note that in this second case, $h_{i,j}(z)$ is necessarily $h_{1,1}(z)$ since the other filters are $0$ at such a
$z$, and we were assuming that we were in the support of $\widehat{\phi}_i(x)$.  Thus, in the second case, we obtain a factor of $|p(z)|$ with $\|z-\frac 3{14}\|<\frac{\delta}4$, and we reset $\delta$ to $\frac1{28}$.  Since $x\in A_n\cup B_n\cup C_n$ satisfies $\frac{8^n}2<|x|<8^{n+1}$, it will take between $n$ and $n+1$ applications of this process of moving from $y$ to $\frac y8$ to move into the unit interval.  From this, if $a_j=1$ for all $j$, we will obtain
$n$ factors of the form $p(z)$ for $\|z-\frac 3{14}\|<\frac 1{112}$.  Thus, by applying the Mean Value Theorem to condition (5) in the definition of the function $p$, we see that we will have
$n$ factors $|h_{i,j}(\frac x{2^l})|$ each of which is less than $(\frac 1{112})^{r+2}$  Each time $a_j=+1$, we will decrease the value of $\delta$ by a factor of $\frac 18$ instead of obtaining a new small factor.  If $a_j=+1$ for all $j$, we will eventually arrive at $h_{1,1}(z)$ for $|z-\frac17|<\frac 1{8^n}$.
Again, by applying the Mean Value Theorem to condition (5) in the definition of $p$, we see that this gives us a factor of $|h_{i,j}(\frac x{2^l})|$ which is less than $(\frac 1{8^n})^{r+2}$.  Values of $x\in A_n\cup B_n\cup C_n$ with a mixture of $a_j=+1$ and $a_j=-1$ will be somewhere between these two extremes, but all will have one or more factors of $p(x)$  among the $|h_{i,j}(\frac x{2^l})|$ with $x$ so close to $\frac 3{14}$ or $\frac17$ that the product of these factors is less than $(\frac 1{8^n})^{r+2}$.  Since the other factors are bounded by $1$, we thus get that each term in
$\widehat{\phi}_i(x)$ is bounded in absolute value by $\frac 1{\sqrt 2}^{3n}(\frac 1{8^n})^{r+2}$.  Thus, for $x\in A_n\cup B_n\cup C_n$, we have
$|\widehat{\phi}_i(x)|<\frac{3^{n+1}}{(\sqrt2)^{3n}8^{n(r+2)}}.$  Finally, for such $x$, $|x^{r+1}\widehat{\phi}_i(x)|<\frac{3^{n+1}8^{(n+1)(r+1)}}{8^{n(r+2)+\frac n2}}
=\frac{3^{n+1}8^{r+1}}{8^\frac{3n}2}\rightarrow 0$ as $n\rightarrow\infty$.
\end{proof}

We are now ready to prove our main theorem:

\begin{theorem}\label{main}  For any fixed integer $r\geq 0$, there exists a real-valued normalized tight frame wavelet $\psi$ such that $\psi\in C^r$, $\widehat{\psi}
\in C^{\infty}$ and such that the multiplicity function associated to $\psi$ is that of the Journ\'e wavelet.
\end{theorem}

\begin{proof}
We build $\psi$ using the technique described in Theorem \ref{bjmp}, with low-pass filter functions $h_{i,j}$ defined in equations (\ref{h11})--(\ref{h22}),
and with high-pass filters given by

$$g_1(x)=\begin{cases}
e^{2\pi i x}\overline{p(x+\frac 12)}&x\in[-\frac 2{7},\frac 2{7})\cr
0&\text{\rm \ otherwise}\end{cases}
$$

$$g_2(x)=\begin{cases}-e^{2\pi i x}\overline{p(x)}&x\in [\frac {-1}7,\frac 1{7})\cr
0&\text{\rm \ otherwise}\end{cases}
$$

It is easily checked that these generalized high-pass filter functions satisfy the filter equations
(\ref{gen4}) and (\ref{gen5}).  In fact, as their form suggests, they were built from the generalized
low-pass filters using the classical high-pass filter obtained from $p$ via equation (\ref{classicalg}).
The frame wavelet $\psi$ is defined as in Theorem \ref{bjmp}, by
 $\widehat{\psi}(x )
\equiv\frac 1{\sqrt{2}}\left(g_{1}(\frac x2 )\widehat{ \phi_1}(\frac x2 )+g_{2}(\frac x2 )\widehat{ \phi_2}(\frac x2 )\right)$,
where  $\widehat{ \phi_1}$ and $\widehat{\phi_2}$ are the functions built from the $\{h_{i,j}\}$ in
Proposition \ref{GMRAex}.  The wavelet $\psi$ is real-valued since the fact that both the low and high pass
filters are symmetric forces $\widehat{\psi}$ to be symmetric as well.  We have that $\widehat{\psi}$
is $C^{\infty}$ since Lemma \ref{Cinfty} established that $\widehat{\phi}_1$ and $\widehat{\phi}_2$
are $C^{\infty}$, and $g_1$ and $g_2$ are as well.  Further, we have
that $|x|^{r+1}\widehat{\psi}(x)\rightarrow 0$ as $x\rightarrow\infty,$ since
$\widehat{\phi}_1$ and $\widehat{\phi}_2$ were shown to have this property in Proposition \ref{xr}
and $g_1$ and $g_2$ are bounded by $\sqrt{2}$.  Thus the wavelet $\psi$ itself satisfies
$\psi\in C^r$.

It remains to show that $\psi$ determines the same GMRA as $\phi_1$ and $\phi_2$, so that this wavelet
also has the multiplicity function of Journ\'e.  It was shown in \cite{BJMP} that $\psi\in V_1$, where
$\{V_j\}$ is the GMRA determined by $\phi_1$ and $\phi_2$ as in Proposition \ref{GMRAex}.  This says that
$\psi$ is {\it obtained from} the GMRA in the sense defined by Zalik \cite{Za}.  We will now show that in
this particular case, $\psi$ is {\it associated with} the GRMA $\{V_j\}$ in the sense
that $V_j$ is the closed linear span of $\{\psi_{l,k}\}_{l<j}$.  (For more on this distinction, see e.g. \cite{Bo}.)
Since we already know
(by $\psi\in V_1$) that $V_j$ contains this closed linear span, it will suffice to show that $\widehat{V_0}$
is contained in the closure of the span of $\{\widehat{\psi}_{l,k}\}_{l<0}$.

We will do this by defining a bounded surjective linear map $T:L^2(S_1)\oplus L^2(S_2)\mapsto\widehat{V_0}$,
that takes a tight frame for $L^2(S_1)\oplus L^2(S_2)$ to $\{\widehat{\psi}_{l,k}\}_{l<0}$.
We define $T(f_1,f_2)=f_1\widehat{\phi_1}+f_2\widehat{\phi_2}$, where $f_i\in L^2(S_i)$ is extended to
a periodic function.  $T$ maps onto
$\widehat{V_0}$ since by Lemma \ref{partialproduct}, the support of $\widehat{\phi_i}$ is the periodization of $S_i$.
To see that $T$ is bounded, we compute the norm

\begin{eqnarray*}
\|T(\{f_i\})\|^2&=& \int_{\mathbb R}\left|f_1(x)\widehat{\phi_1}(x)+f_2(x)\widehat{\phi_2}(x)\right|^2 dx\\
&=& \int_{\mathbb R}\left|f_1(x)\widehat{\phi_1}(x)\right|^2+\left|f_2(x)\widehat{\phi_2}(x)\right|^2 dx\\
&=&\sum_{j\in\mathbb Z}\int_{\mathbb T} \left|f_1(x)\widehat{\phi_1}(x+j)\right|^2+\left|f_2(x)\widehat{\phi_2}(x+j)\right|^2dx\\
&=&\int_{\mathbb T}\left|f_1(x)\right|^2\sum_{j\in\mathbb Z}\left|\widehat{\phi_1}(x+j)\right|^2 +\left|f_2(x)\right|^2\sum_{j\in\mathbb Z}\left|\widehat{\phi_2}(x+j)\right|^2dx,
\end{eqnarray*}
where the second equality follows from the disjointness of the supports of $\widehat{\phi_1}$ and $\widehat{\phi_2}$.
Let $\textup{Per}\, \phi_i(x)\equiv\sum_{j\in\mathbb Z}\left|\widehat{\phi_i}(x+j)\right|^2.$
To show that $T$ is bounded, it will suffice to show that $\textup{Per}\,\phi_i$ is bounded on $\mathbb T$ for $i=1,2$.  This
follows from  Lemma \ref{partialproduct} by writing, for $x\in[-\frac12,\frac12]$,
\begin{eqnarray*}
\textup{Per}\, \phi_i(x)&=&\sum_{n=0}^{\infty}\sum_{\{a_j\}\in\{-1,1\}^{n}}\widehat{\phi}_i(x + \sum_{j=0}^{n}a_j 8^j)\\
&\leq&\sum_{n=0}^{\infty}2^n\frac1{8^{(n-1)r}}
\end{eqnarray*}
where the last step follows from Proposition \ref{xr}.

Now to see that $T$ takes a tight frame for $L^2(S_1)\oplus L^2(S_2)$ to $\{\widehat{\psi}_{l,k}\}_{l<0}$,
we adopt some notation from \cite{BJMP}.  Define the operator $S_G:L^2(\mathbb{T})\rightarrow L^2(S_1)\oplus L^2(S_2)$ by
 $$[S_G(f)](x)=\{g_1(x)f(2x),g_2(x)f(2x)\}$$
and
$S_H:L^2(S_1)\oplus L^2(S_2)\rightarrow L^2(S_1)\oplus L^2(S_2)$ by
$$[S_H(\{f_1,f_2\})](x) = \{h_{1,1}(x)f_1(2x)+h_{2,1}(x)f_2(2x), h_{1,2}(x)f_1(2x)\}.$$
It is shown in \cite{BJMP} (Lemma 3.3) that $\{S_H^n(S_G(e^{2\pi i k x}))\}_{k\in\mathbb{Z}, n=0,1,2,\cdots}$ forms a normalized tight frame
for $L^2(S_1)\oplus L^2(S_2)$.  We complete the proof by showing that the operator $T$ takes this frame to $\{\widehat{\psi}_{l,k}\}_{l<0}.$
Note first that for $n=0$, we have
\begin{eqnarray*}
T(S_G(e^{2\pi i k x}))&=&g_1(x)e^{4\pi i k x}\widehat{\phi_1}(x)+g_2(x)e^{4\pi i k x}\widehat{\phi_1}(x)\\
&=&e^{4\pi i k x}\widehat{\psi}(2x)\\
&=&\widehat{\psi}_{-1,k}
\end{eqnarray*}
Thus the result will follow from
the fact that $T$ intertwines the operators $S_H$ on $L^2(S_1)\oplus L^2(S_2)$ and inverse dilation $\delta^{-1}$ on $\widehat{V_0}$:

\begin{eqnarray*} [T(S_H(\{f_1,f_2\}))](x)
&=& h_{1,1}(x)f_1(2x)+h_{2,1}(x)f_2(2x)\widehat{\phi_1}(x)+ h_{1,2}(x)f_1(2x)\widehat{\phi_2}(x)\cr
&=&  f_1(2x)\left( h_{1,1}(x)\widehat{\phi_1}(x)+h_{1,2}(x)\widehat{\phi_2}(x)\right)+ f_2(2x)h_{2,1}(x)\widehat{\phi_1}(x) \cr
&=& f_1(2x) \sqrt 2 \widehat {\phi_1}(2x)+f_2(2x) \sqrt 2 \widehat {\phi_2}(2x) \cr
&=& \left[\delta^{-1}\left( f_1\widehat{\phi_1} +f_2\widehat{\phi_2}\right)\right](x)\cr
&=& [\delta^{-1}(T(f))](x).\end{eqnarray*}
\end{proof}

\begin{remark}
While the dilates and translates of the wavelet $\psi$ defined in Theorem \ref{main} do form a normalized tight frame for all of $L^2(\mathbb{R})$,
the translates of  $\psi$ do not even form a frame for their linear span.  We can see this using the result in \cite{BLi} once again, since
$\textup{Per}\, \psi(x)\equiv\sum_{j\in\mathbb Z}\left|\widehat{\psi}(x+j)\right|^2$ is not bounded away from $0$ on its support.  In
particular, $\textup{Per}\, \psi$ is continuous, yet nonzero on $\pm(\frac 27, \frac 47)$, and 0 on $\pm(\frac 47,\frac 67)$.
\end{remark}

\end{document}